\input amstex
\input amsppt.sty   
\hsize 11cm
\vsize 16.6cm
\magnification=\magstep1
\NoBlackBoxes

\define\nmb#1#2{#2}      
\def\ign#1{}             
\redefine\o{\circ}
\define\X{\frak X}
\define\al{\alpha}

\define\ga{\gamma}
\define\de{\delta}

\define\la{\lambda}
\define\rh{\rho}
\define\si{\sigma}

\define\Om{\Omega}
\redefine\i{^{-1}}
\define\row#1#2#3{#1_{#2},\ldots,#1_{#3}}
\define\x{\times}
\redefine\L{{\Cal L}}
\define\g{{\frak g}}
\define\ome{\varOmega}           
\define\grad{\operatorname{grad}^\ome }
\define\a{\omega}                 
\define\Lip{\operatorname{Lip}}

\def\today{\ifcase\month\or
 January\or February\or March\or April\or May\or June\or
 July\or August\or September\or October\or November\or December\fi
 \space\number\day, \number\year}
\topmatter
\title All Unitary Representations Admit Moment Mappings \endtitle
\author  Peter W. Michor  \endauthor
\address{Institut f\"ur Mathematik, Universit\"at Wien,
Strudlhofgasse 4, A-1090 Wien, Austria.}\endaddress
\endtopmatter
\document

\heading\nmb0{1}. Calculus of smooth mappings \endheading

\subheading{\nmb.{1.1}} The traditional differential calculus works 
well for finite dimensional vector spaces and for Banach spaces. For 
more general locally convex spaces a whole flock of different 
theories were developed, each of them rather complicated and none 
really convincing. The main difficulty is that the composition of 
linear mappings stops to be jointly continuous at the level of Banach 
spaces, for any compatible topology. This was the original motivation 
for the development of a whole new field within general topology, 
convergence spaces.

Then in 1982, Alfred Fr\"olicher and Andreas Kriegl presented 
independently the solution to the question for the right differential 
calculus in infinite dimensions. They joined forces in the further 
development of the theory and the (up to now) final outcome is the 
book \cite{F-K}.

In this section I will sketch the basic definitions and the most 
important results of the Fr\"olicher-Kriegl calculus.

\subheading{\nmb.{1.2}. The $c^\infty$-topology} Let $E$ be a 
locally convex vector space. A curve $c:\Bbb R\to E$ is called 
{\it smooth} or $C^\infty$ if all derivatives exist and are 
continuous - this is a concept without problems. Let 
$C^\infty(\Bbb R,E)$ be the space of smooth functions. It can be 
shown that $C^\infty(\Bbb R,E)$ does not depend on the locally convex 
topology of $E$, only on its associated bornology (system of bounded 
sets).

The final topologies with respect to the following sets of mappings 
into E coincide:
\roster
\item $C^\infty(\Bbb R,E)$.
\item Lipschitz curves (so that $\{\frac{c(t)-c(s)}{t-s}:t\neq s\}$ 
     is bounded in $E$). 
\item $\{E_B\to E: B\text{ bounded absolutely convex in }E\}$, where 
     $E_B$ is the linear span of $B$ equipped with the Minkowski 
     functional $p_B(x):= \inf\{\la>0:x\in\la B\}$.
\item Mackey-convergent sequences $x_n\to x$ (there exists a sequence 
     $0<\la_n\nearrow\infty$ with $\la_n(x_n-x)$ bounded).
\endroster
This topology is called the $c^\infty$-topology on $E$ and we write 
$c^\infty E$ for the resulting topological space. In general 
(on the space $\Cal D$ of test functions for example) it is finer 
than the given locally convex topology, it is not a vector space 
topology, since scalar multiplication is no longer jointly 
continuous. The finest among all locally convex topologies on $E$ 
which are coarser than $c^\infty E$ is the bornologification of the 
given locally convex topology. If $E$ is a Fr\'echet space, then 
$c^\infty E = E$. 

\subheading{\nmb.{1.3}. Convenient vector spaces} Let $E$ be a 
locally convex vector space. $E$ is said to be a {\it convenient 
vector space} if one of the following equivalent (completeness) 
conditions is satisfied:
\roster
\item Any Mackey-Cauchy-sequence (so that $(x_n-x_m)$ is Mackey 
     convergent to 0) converges. This is also called 
     $c^\infty$-complete.
\item If $B$ is bounded closed absolutely convex, then $E_B$ is a 
     Banach space.
\item Any Lipschitz curve in $E$ is locally Riemann integrable.
\item For any $c_1\in C^\infty(\Bbb R,E)$ there is 
     $c_2\in C^\infty(\Bbb R,E)$ with $c_1'=c_2$ (existence of 
     antiderivative).
\endroster

\proclaim{\nmb.{1.4}. Lemma} Let $E$ be a locally convex space.
Then the following properties are equivalent:
\roster
\item $E$ is $c^\infty$-complete.
\item If $f:\Bbb R^k\to E$ is scalarwise $\Lip^k$, then $f$ is 
     $\Lip^k$, for $k>1$.
\item If $f:\Bbb R\to E$ is scalarwise $C^\infty$ then $f$ is 
     differentiable at 0.
\item If $f:\Bbb R\to E$ is scalarwise $C^\infty$ then $f$ is 
     $C^\infty$.
\endroster
\endproclaim
Here a mapping $f:\Bbb R^k\to E$ is called $\Lip^k$ if all partial 
derivatives up to order $k$ exist and are Lipschitz, locally on 
$\Bbb R^n$. $f$ scalarwise $C^\infty$ means that $\la\o f$ is $C^\infty$  
for all continuous linear functionals on $E$.

This lemma says that a convenient vector space one can recognize 
smooth curves by investigating compositions with continuous linear 
functionals.

\subheading{\nmb.{1.5}. Smooth mappings} Let $E$ and $F$ be locally 
convex vector spaces. A mapping $f:E\to F$ is called {\it smooth} or 
$C^\infty$, if $f\o c\in C^\infty(\Bbb R,F)$ for all 
$c\in C^\infty(\Bbb R,E)$; so 
$f_*: C^\infty(\Bbb R,E)\to C^\infty(\Bbb R,F)$ makes sense.
Let $C^\infty(E,F)$ denote the space of all smooth mapping from $E$ 
to $F$.

For $E$ and $F$ finite dimensional this gives the usual notion of 
smooth mappings: this has been first proved in \cite{Bo}.
Constant mappings are smooth. Multilinear mappings are smooth if and 
only if they are boun\-ded. Therefore we denote by $L(E,F)$ the space 
of all bounded linear mappings from $E$ to $F$.

\subheading{\nmb.{1.6}. Structure on $C^\infty(E,F)$} We equip the 
space $C^\infty(\Bbb R,E)$ with the bornologification of the topology 
of uniform convergence on compact sets, in all derivatives 
separately. Then we equip the space $C^\infty(E,F)$ with the 
bornologification of the initial topology with respect to all 
mappings $c^*:C^\infty(E,F)\to C^\infty(\Bbb R,F)$, $c^*(f):=f\o c$, 
for all $c\in C^\infty(\Bbb R,E)$.

\proclaim{\nmb.{1.7}. Lemma } For locally convex spaces $E$ and $F$ 
we have:
\roster
\item If $F$ is convenient, then also $C^\infty(E,F)$ is convenient, 
     for any $E$. The space $L(E,F)$ is a closed linear subspace of 
     $C^\infty(E,F)$, so it also  convenient.
\item If $E$ is convenient, then a curve $c:\Bbb R\to L(E,F)$ is 
     smooth if and only if $t\mapsto c(t)(x)$ is a smooth curve in $F$ 
     for all $x\in E$.
\endroster
\endproclaim

\proclaim{\nmb.{1.8}. Theorem} The category of convenient vector 
spaces and smooth mappings is cartesian closed. So we have a natural 
bijection 
$$C^\infty(E\x F,G)\cong C^\infty(E,C^\infty(F,G)),$$
which is even a diffeomorphism.
\endproclaim

Of coarse this statement is also true for $c^\infty$-open subsets of 
convenient vector spaces. 

\proclaim{\nmb.{1.9}. Corollary } Let all spaces be convenient vector 
spaces. Then the following canonical mappings are smooth.
$$\align
&\operatorname{ev}: C^\infty(E,F)\x E\to F,\quad 
     \operatorname{ev}(f,x) = f(x)\\
&\operatorname{ins}: E\to C^\infty(F,E\x F),\quad
     \operatorname{ins}(x)(y) = (x,y)\\
&(\quad)^\wedge :C^\infty(E,C^\infty(F,G))\to C^\infty(E\x F,G)\\
&(\quad)\spcheck :C^\infty(E\x F,G)\to C^\infty(E,C^\infty(F,G))\\
&\operatorname{comp}:C^\infty(F,G)\x C^\infty(E,F)\to C^\infty(E,G)\\
&C^\infty(\quad,\quad):C^\infty(F,F')\x C^\infty(E',E)\to 
     C^\infty(C^\infty(E,F),C^\infty(E',F'))\\
&\qquad (f,g)\mapsto(h\mapsto f\o h\o g)\\
&\prod:\prod C^\infty(E_i,F_i)\to C^\infty(\prod E_i,\prod F_i)
\endalign$$
\endproclaim

\proclaim{\nmb.{1.10}. Theorem} Let $E$ and $F$ be convenient vector 
spaces. Then the differential operator 
$$\gather d: C^\infty(E,F)\to C^\infty(E,L(E,F)), \\
df(x)v:=\lim_{t\to0}\frac{f(x+tv)-f(x)}t,
\endgather$$
exists and is linear and bounded (smooth). Also the chain rule holds: 
$$d(f\o g)(x)v = df(g(x))dg(x)v.$$
\endproclaim

\subheading{\nmb.{1.11}. Remarks } Note that the conclusion of 
theorem \nmb!{1.8} is the starting point of the classical calculus of 
variations, where a smooth curve in a space of functions was assumed 
to be just a smooth function in one variable more.

If one wants theorem \nmb!{1.8} to be true and assumes some other obvious 
properties, then the calculus of smooth functions is already uniquely 
determined.

There are, however, smooth mappings which are not continuous. This is 
unavoidable and not so horrible as it might appear at first sight. 
For example the evaluation $E\x E'\to\Bbb R$ is jointly continuous if 
and only if $E$ is normable, but it is always smooth. Clearly smooth 
mappings are continuous for the $c^\infty$-topology.

For Fr\'echet spaces smoothness in the sense described here coincides 
with the notion $C^\infty_c$ of \cite{Ke}. This is the 
differential calculus used by \cite{Mic1}, \cite{Mil}, 
and \cite{P-S}.

A prevalent opinion in contemporary mathematics is,
that for infinite
dimensional calculus each serious application needs its own 
foundation. By a serious application one obviously means some 
application of
a hard inverse function theorem. These theorems can be proved,
if by assuming enough a priori estimates one creates enough
Banach space situation for some modified iteration procedure to
converge. Many authors try to build their platonic idea of an a
priori estimate into their differential calculus. I think that this 
makes the calculus
inapplicable and hides the origin of the a priori estimates. I
believe, that the calculus itself should be as easy to use as
possible, and that all further assumptions (which most often
come from ellipticity of some nonlinear partial differential
equation of geometric origin) should be treated separately, in a
setting depending on the specific problem. I am sure that in
this sense the Fr\"olicher-Kriegl calculus as presented here and its 
holomorphic and real analytic offsprings in sections \nmb!{2} and 
\nmb!{3} below are universally usable for most applications.

Let me point out as a final remark, that also the cartesian closed 
calculus for holomorphic mappings along the same lines is available 
in \cite{K-N}, and recently the cartesian closed calculus for real 
analytic mapping was developed in \cite{K-M}.

\heading \nmb0{2}. The moment mapping for unitary representations 
\endheading

The following is a review of the results 
obtained in \cite{Mic2}. We include only one proof, the central 
application of the Fr\"olicher-Kriegl calculus.

\subheading{\nmb.{2.1}} Let $G$ be any (finite dimensional second
countable) real Lie group, and let $\rh: G\to  U(\bold H)$ be a
unitary representation on a Hilbert space $\bold H$. Then the
associated mapping $\hat \rh: G\x \bold H \to  \bold H$ is in general
{\it not} jointly continuous, it is only separately continuous,
so that $g\mapsto \rh(g)x$, $G\to  \bold H$, is continuous for any
$x\in \bold H$.

\subheading{Definition} A vector $x \in \bold H$ is called {\it
smooth} (or {\it real analytic}) if the mapping $g\mapsto
\rh(g)x$, $G\to \bold H$ is smooth (or real analytic).
Let us denote by $\bold H_\infty$ the linear subspace of all
smooth vectors in $\bold H$. Then we have an embedding 
$j:\bold H_\infty \to  C^\infty(G,\bold H)$, given by 
$x\mapsto(g\mapsto\rh(g)x)$. We equip $C^\infty(G,\bold H)$ with
the compact $C^\infty$-topology (of uniform convergence on
compact subsets of $G$, in all derivatives separately). Then
it is easily seen (and proved in \cite {Wa, p 253}) that 
$\bold H_\infty$ is a closed linear subspace. So with the
induced topology $\bold H_\infty$ becomes a Fr\`echet space.
Clearly $\bold H_\infty$ is also an invariant subspace, so we
have a representation $\rh: G\to  L(\bold H_\infty,\bold H_\infty)$.
For more detailed information on $\bold H_\infty$ see
\cite{Wa, chapt. 4.4.} or \cite{Kn, chapt. III.}.

\proclaim{\nmb.{2.2}. Theorem} The mapping $\hat\rh: G\x\bold H_\infty \to 
\bold H_\infty$ is smooth.
\endproclaim
\demo{Proof} By cartesian closedness of the
Fr\"olicher-Kriegl calculus \nmb!{1.8} it suffices to
show that the canonically associated mapping 
$$\hat\rh\spcheck :G \to  C^\infty(\bold H_\infty,\bold H_\infty)$$
is smooth; but it takes values in the closed subspace 
$L(\bold H_\infty,\bold H_\infty)$ 
of all bounded linear operators. So by  it suffices to show that the mapping 
$\rh:G \to  L(\bold H_\infty,\bold H_\infty)$
is smooth. But for that, since $\bold H_\infty$ is a Fr\`echet
space, thus convenient in the sense of Fr\"olicher-Kriegl, 
by \nmb!{1.7}(2) it suffices to show that 
$$G @>{\rho}>>  L(\bold H_\infty,\bold H_\infty) @>{ev_x}>> \bold H_\infty$$ 
is smooth for each $x \in \bold
H_\infty$. This requirement means that $g\mapsto \rh(g)x$,
$G\to \bold H_\infty$, is smooth.  For this it suffices to show
that $$\align &G\to \bold H_\infty @>j>> C^\infty(G,\bold H),\\
&g\mapsto \rh(g)x \mapsto(h\mapsto \rh(h)(g)x),\endalign$$ is
smooth.  But again by cartesian closedness it suffices to show
that the associated mapping $$\align &G\x G\to  \bold H,\\
&(g,h)\mapsto \rh(h)(g)x = \rh(hg)x,\endalign$$ is smooth. And
this is the case since $x$ is a smooth vector.
\qed\enddemo

\subheading{\nmb.{2.3}} we now consider $\bold H_\infty$ as a "weak"
symplectic Fr\`echet manifold, equipped with the symplectic
structure $\ome$, the restriction of the imaginary part of the
Hermitian inner product $\langle \quad,\quad\rangle$ on $\bold H$. Then
$\ome\in \Om^2(\bold H_\infty)$ is a closed 2-form which is non
degenerate in the sense that 
$$\check \ome: T\bold H_\infty = \bold H_\infty\x\bold H_\infty \to 
T^*\bold H_\infty = \bold H_\infty\x{\bold H_\infty}'$$
is injective (but not surjective), where 
$\bold H_\infty{}' = L(\bold H_\infty,\Bbb R)$
denotes the real topological dual space. This is the meaning of "weak"
above.

\subheading{\nmb.{2.4}. Review} For a finite dimensional symplectic
manifold $(M,\ome)$ we have the following exact sequence of
Lie algebras:
$$0\to  H^0(M)\to C^\infty(M)  @>{\grad}>> \X_\ome(M) @>\ga>> H^1(M)\to 0$$
Here $H^*(M)$ is the real De Rham cohomology of $M$, the space
$C^\infty(M)$ is equipped with the Poisson bracket $\{\quad,\quad\}$,
$\X_\ome(M)$ consists of all vector fields $\xi$ with $\L_\xi\ome=0$
(the locally Hamiltonian vector fields), which is a Lie algebra
for the Lie bracket. 
$\grad f$ is the Hamiltonian vector field
for $f\in C^\infty(M)$ given by $i(\grad f)\ome = df$, and 
$\ga(\xi) = [i_\xi\ome]$. The spaces $H^0(M)$ and $H^1(M)$ are
equipped with the zero bracket.

Given a symplectic left action $\ell:G\x M\to  M$ of a connected
Lie group $G$ on $M$, the first partial derivative of $\ell$
gives a mapping $\ell':\g \to  \X_\ome(M)$ which sends each element
$X$ of the Lie algebra $\g$ of $G$ to the fundamental vector
field. This is a Lie algebra homomorphism. 
$$\CD
H^0(M) @>i>> C^\infty(M)  @>{\grad}>> \X_\ome(M) @>\ga>> H^1(M) \\
@.           @A\si AA            @AA{\ell'}A      @.  \\
       @.    \g           @=    \g        @.            
\endCD$$
A linear lift $\si:\g\to  C^\infty(M)$ of $\ell'$  with 
$\grad\o\si=\ell'$ exists if
and only if $\ga\o\ell'=0$ in $H^1(M)$. This lift $\si$ may be
changed to a Lie algebra homomorphism if and only if the
$2$-cocycle $\bar\si:\g\x\g\to  H^0(M)$, given by 
$(i\o\bar\si)(X,Y) = \{\si(X),\si(Y)\} - \si([X,Y])$, vanishes
in $H^2(\g,H^0(M))$, for if $\bar\si = \de\al$ then $\si-i\o\al$
is a Lie algebra homomorphism.

If $\si:\g\to C^\infty(M)$ is a Lie algebra homomorphism, we may
associate the {\it moment mapping} 
$\mu:M\to \g'=L(\g,\Bbb R)$ to it, which is 
given by
$\mu(x)(X) = \si(X)(x)$ for $x\in M$ and $X\in \g$.
It is $G$-equivariant for a suitably chosen (in general affine)
action of $G$ on $\g'$.
See \cite{We} or \cite{L-M} for all this.

\subheading{\nmb.{2.5}} We now want to carry over to the setting of \nmb!{2.1}
and \nmb!{2.2} the procedure of \nmb!{2.4}. The first thing to note is that
the hamiltonian mapping 
$\grad:C^\infty(\bold H_\infty) \to  \X_\ome(\bold H_\infty)$
does not make sense in general, since 
$\check\ome:\bold H_\infty \to  {\bold H_\infty}'$ is not invertible:
$\grad f = \check\ome\i df$ is defined only for those $f\in
C^\infty(\bold H_\infty)$ with $df(x)$ in the image of $\check\ome$ for all
$x\in \bold H_\infty$. A similar difficulty arises for the
definition of the Poisson bracket on $C^\infty(\bold H_\infty)$.

Let $\langle x,y\rangle = Re\langle x,y\rangle +
\sqrt{-1}\ome(x,y)$ be the decomposition of the hermitian inner
product into real and imaginary parts. Then $Re\langle
x,y\rangle = \ome(\sqrt{-1}x,y)$, thus the real linear subspaces
$\check\ome(\bold H_\infty) = \ome(\bold H_\infty,\quad)$ and 
$Re\langle \bold H_\infty,\quad\rangle$ of $\bold H_\infty{}' =
L(\bold H_\infty,\Bbb R)$ coincide.

\subheading{\nmb.{2.6} Definition} Let $\bold H_\infty^*$ denote
the real linear subspace 
$$\bold H_\infty^*=\ome(\bold H_\infty,\quad) = 
Re\langle \bold H_\infty,\quad\rangle$$ 
of $\bold H_\infty{}' = 
L(\bold H_\infty,\Bbb R)$, and let us call it the {\it smooth dual}
of $\bold H_\infty$ in view of the embedding of test functions
into distributions. We have two canonical isomorphisms 
$\bold H_\infty^* \cong \bold H_\infty$ induced by $\ome$ and 
$Re\langle\quad,\quad\rangle$, respectively. Both induce the
same Fr\'echet topology on $\bold H_\infty^*$, which we fix from
now on.

\subheading{\nmb.{2.7} Definition} Let 
$C^\infty_*(\bold H_\infty,\Bbb R) \subset
C^\infty(\bold H_\infty,\Bbb R)$ 
denote the linear subspace consisting of all smooth functions
$f:\bold H_\infty \to  \Bbb R$ such that each iterated derivative 
$d^kf(x)\in L^k_{\text{sym}}(\bold H_\infty,\Bbb R)$ has the
property that 
$$d^kf(x)(\quad,\row y2k) \in {\bold H_\infty}^*$$
is actually in the smooth dual $\bold H_\infty^* \subset
{\bold H_\infty}'$ for all $x,\row y2k \in \bold H_\infty$, and
that the mapping 
$$\gather \prod^k\bold H_\infty \to  \bold H_\infty \\
(x,\row y2k)\mapsto \check\ome\i(df(x)(\quad,\row y2k))
\endgather$$
is smooth. Note that we could also have used
$Re\langle\quad,\quad\rangle$ instead of $\ome$. By the symmetry
of higher derivatives this is then true for all entries of
$d^kf(x)$, for all $x$.  

\proclaim{\nmb.{2.8} Lemma} For $f \in C^\infty(\bold H_\infty,\Bbb R)$ the
following assertions are equivalent:
\roster
\item $df:\bold H_\infty \to  \bold H_\infty{}'$ factors to a
        smooth mapping $\bold H_\infty \to  \bold H_\infty^*$.
\item $f$ has a smooth $\ome$-gradient 
        $\grad f\in \X(\bold H_\infty) = C^\infty(\bold
        H_\infty,\bold H_\infty)$ such that 
        $df(x)y = \ome(\grad f(x),y)$.
\item $f \in C^\infty_*(\bold H_\infty,\Bbb R)$.
\endroster
\endproclaim

\proclaim{\nmb.{2.9}. Theorem} 
The mapping 
$$\grad :C^\infty_*(\bold H_\infty,\Bbb R) 
\to  \X_\ome(\bold H_\infty),\qquad \grad f := \check\ome\i\o df,$$
is well defined; also the
Poisson bracket 
$$\gather
\{\quad,\quad\}: C^\infty_*(\bold H_\infty,\Bbb R) \x
C^\infty_*(\bold H_\infty,\Bbb R) \to  
C^\infty_*(\bold H_\infty,\Bbb R),\\
\{f,g\}:= i(\grad f)i(\grad g)\ome = \ome(\grad g,\grad f) =\\
= (\grad f)(g) = dg(\grad f) 
\endgather$$ 
is well defined and gives a Lie algebra structure to
the space $C^\infty_*(\bold H_\infty,\Bbb R)$.

We also have the following long exact sequence of Lie algebras
and Lie algebra homomorphisms:
$$0\to  H^0(\bold H_\infty)\to C^\infty_*(\bold H_\infty,\Bbb R)  
@>{\grad}>> \X_\ome(\bold H_\infty) @>\ga>> H^1(\bold H_\infty) = 0$$ 
\endproclaim

\subheading{\nmb.{2.10}} We consider now again as in
\nmb!{2.1} a unitary representation $\rh: G\to  U(\bold H)$. By theorem
\nmb!{2.2} the associated mapping 
$\hat \rh: G\x \bold H_\infty \to  \bold H_\infty$ 
is smooth, so we have the infinitesimal mapping 
$\rh':\g \to  \X(\bold H_\infty)$, given by 
$\rh'(X)(x) = T_e(\hat\rho(\quad,x))$ for $X\in \g$ and $x\in
\bold H_\infty$. Since $\rho$ is a unitary
representation, the mapping $\rho'$ has values in the Lie
subalgebra of all linear hamiltonian vector fields $\xi \in
\X(\bold H_\infty)$ which respect the symplectic form $\ome$,
i.e. $\xi:\bold H_\infty \to  \bold H_\infty$ is linear and
$\L_\xi\ome = 0$.  

Now let us consider the mapping 
$\check\ome\o \rho'(X): \bold H_\infty \to  T(\bold H_\infty) \to 
T^*(\bold H_\infty)$. 
We have $d(\check\ome\o\rho'(X)) = d(i_{\rho'(X)}\ome) =
\L_{\rho'(X)}\ome = 0$, so the linear 1-form
$\check\ome\o\rho'(X)$ is closed, and since $H^1(\bold H_\infty)
= 0$, it is exact. So there is a function 
$\si(X)\in C^\infty(\bold H_\infty,\Bbb R)$ with $d\si(X) =
\check\ome\o\rho'(X)$, 
and $\si(X)$ is uniquely determined up to addition of a constant.
If we require $\si(X)(0) = 0$, then $\si(X)$ is uniquely
determined and is a quadratic function. In fact we have 
$\si(X)(x) = \int_{c_x}\check\ome\o\rho'(X)$, where $c_x(t) =
tx$. Thus
$$\align 
\si(X)(x) 
&= {\tsize\int_0^1} \ome(\rho'(X)(tx),\tfrac d{dt}tx)dt = \\
&= \ome(\rho'(X)(x),x){\tsize\int_0^1dt} \\
&= \tfrac 12 \ome(\rho'(X)(x),x).
\endalign$$

\proclaim{\nmb.{2.11}. Lemma} 
The mapping 
$$\si:\g \to  C^\infty_*(\bold H_\infty,\Bbb R),\qquad
\si(X)(x) = \tfrac12\ome(\rho'(X)(x),x)$$ 
for $X\in \g$ and $x\in
\bold H_\infty$, is a Lie algebra homomorphism
and $\grad\o\si = \rho'$.

For $g\in G$ we have $\rho(g)^*\si(X)=\si(X)\o \rho(g) =
\si(Ad(g\i)X)$, so $\si$ is $G$-equivariant.
\endproclaim

\subheading{\nmb.{2.12}. The moment mapping} For a unitary
representation $\rho: G \to  U(\bold H)$ we can now define the
{\it moment mapping} 
$$\gather 
\mu: \bold H_\infty \to  \g' = L(\g,\Bbb R), \\
\mu(x)(X) := \si(X)(x) = \tfrac12 \ome(\rho'(X)x,x),
\endgather$$
for $x \in \bold H_\infty$ and $X \in \g$.

\proclaim{\nmb.{2.13} Theorem} The moment mapping 
$\mu:\bold H_\infty \to  \g'$ has the following properties:
\roster
\item $(d\mu(x)y)(X)=\ome(\rho'(X)x,y)$ for 
    $x,y\in \bold H_\infty$ and $X\in\g$, so 
    $\mu\in C^\infty_*(\bold H_\infty,\g')$.
\item For $x\in\bold H_\infty$ the image of 
    $d\mu(x):\bold H_\infty\to\g'$ is the
    annihilator $\g_x^\ome$ of the Lie algebra 
    $\g_x = \{X\in\g:\rho'(X)(x)=0\}$  of the isotropy group
    $G_x = \{g\in G:\rho(g)x=x\}$ in $\g'$.
\item For $x\in\bold H_\infty$ the kernel of $d\mu(x)$ is
    $$(T_x(\rho(G)x))^\ome = 
    \{y\in\bold H_\infty:\ome(y,T_x(\rho(G)x))=0\},$$
    the $\ome$-annihilator
    of the tangent space at $x$ of the $G$-orbit through $x$.
\item The moment mapping is equivariant:
    $Ad'(g)\o \mu = \mu \o \rho(g)$ for all $g\in G$,
    where $Ad'(g)=Ad(g\i)':\g'\to\g'$ is the coadjoint action.
\item The pullback operator 
    $\mu^*: C^\infty(\g,\Bbb R) \to C^\infty(\bold H_\infty,\Bbb R)$ 
    actually has values in the subspace 
    $C^\infty_*(\bold H_\infty,\Bbb R)$.
    It also is a Lie algebra homomorphism for the Poisson
    brackets involved.
\endroster
\endproclaim

\subheading{\nmb.{2.14}} Let again $\rho: G\to U(\bold H)$ be a
unitary representation of a Lie group $G$ on a Hilbert space
$\bold H$. 

\demo{Definition} A vector $x\in \bold H$ is called {it real
analytic} if the mapping $g\mapsto \rho(g)x$, $G\to \bold H$ is
a real analytic mapping, in the real analytic structure of the
Lie group $G$. 
\enddemo

We will use from now on the theory of real analytic mappings in
infinite dimensions as developed in \cite{K-M}.
So the following conditions on $x\in\bold H$ are equivalent:
\roster
\item $x$ is a real analytic vector.
\item $\g\ni X\mapsto \rho(\exp X)x$ is locally near 0
    given by a converging power series.
\item For each $y\in \bold H$ the mapping $\g\ni X\mapsto
    \langle \rho(\exp X)x,y\rangle\in \Bbb C$ is smooth and real
    analytic along affine lines in $\g$, locally near 0.
\endroster
The only nontrivial part is \therosteritem3 $\Rightarrow$
\therosteritem1, and this follows from \cite{K-M, 1.6 and 2.7} 
and the fact, that $\rho$ is a representation.

Let $\bold H_\a$ denote the vector space of all
real analytic vectors in $\bold H$. 
Then we have a linear
embedding $j:\bold H_\a\to C^\a(G,\bold H)$ into the space of real
analytic mappings, given by 
$x\mapsto (g\mapsto \rho(g)x)$. We equip $C^\a(G,\bold H)$ with
the convenient vector space structure described in \cite{K-M,
5.4, see also 3.13}. 
Then $\bold H_\a$ consists of all equivariant functions in
$C^\a(G,\bold H)$ and is therefore a closed subspace. So it is
a convenient vector space with the induced structure.

The space $\bold H_\a$ is dense in the Hilbert space $\bold H$
by \cite{Wa, 4.4.5.7} and an invariant subspace, so we have a
representation $\rho: G\to L(\bold H_\a,\bold H_\a)$.

\proclaim{\nmb.{2.15}. Theorem} The mapping 
$\hat\rho:G\x\bold H_\a\to \bold H_\a$ is real analytic in the
sense of \rm \cite {K-M}.
\endproclaim
\demo{Proof} Similar to the proof of theorem \nmb!{2.2}.
\qed\enddemo

\subheading{\nmb.{2.16}} Again we consider now $\bold H_\a$ as a
"weak" symplectic real analytic Fr\'echet manifold, equipped
with the symplectic structure $\ome$, the restriction of the
imaginary part of the hermitian inner product $\langle
\quad,\quad\rangle$ on $\bold H$. Then again $\ome\in \Om^2(\bold
H_\a)$  is a closed 2-form which is non degenerate in the sense that 
$\check\ome:\bold H_\a\to \bold H_\a'=L(\bold H_\a,\Bbb R)$ 
is injective. Let 
$$\bold H_\a^*:= \check\ome(\bold H_\a)=\ome(\bold H_\a,\quad)=
Re\langle \bold H_\a,\quad\rangle\subset \bold H_\a'=L(\bold
H_\a,\Bbb R)$$
again denote the {\it analytic dual} of $\bold H_\a$, equipped
with the topology induced by the isomorphism
with $\bold H_\a$.

\subheading{\nmb.{2.17} Remark} All the results leading to the
smooth moment mapping can now be carried over to the real
analytic setting with {\it no} changes in the proofs. So all
statements from \nmb!{2.9} to \nmb!{2.13} are valid in the real
analytic situation. We summarize this in one more result:

\proclaim{\nmb.{2.18} Theorem} Consider the injective linear
continuous $G$-equivariant mapping $i:\bold H_\a\to \bold
H_\infty$. Then for the smooth moment mapping 
$\mu:\bold H_\infty\to \g'$ from \nmb!{2.13} the composition
$\mu\o i:\bold H_\a\to \bold H_\infty\to\g'$ is real analytic.
It is called the real analytic moment mapping.
\endproclaim

\subheading{\nmb.{2.19}. Remarks}
It is my conjecture that for an irreducible representation which
is constructed by geometric quantization of an coadjoint orbit
(the Kirillov method), the restriction of the moment mapping to
the intersection of the unit sphere with the space of smooth
vectors takes values exactly in the orbit one started with, if
the construction is suitably normalized. 

I have been unable to
prove this conjecture in general, but Herbert Wiklicky \cite{Wi}
has checked that this is true for the
Heisenberg group. He also checked that this moment mapping
produces the expectation value for the (angular) momentum in
physically relevant situations and he claims that this moment
mapping describes a sort of classical limit for the quantum
theory described by the unitary representation in question.

Let me add some thoughts on the r\^ole of the moment mapping in
the study of unitary representations. I think that its
restriction to the intersection of the unit sphere with the
space of smooth vectors maps to one coadjoint orbit, if the
representation is irreducible (I was unable to prove this). It
is known that not all irreducible representations come from line
bundles over coadjoint orbits (alias geometric quantization),
but there might be a higher dimensional vector bundle over this
coadjoint orbit, whose space of sections contains the space of
smooth vectors as subspace of sections which are covariantly
constant along some complex polarization.

\enddocument